\documentclass[a4paper,twoside,11pt]{article}

\usepackage[USenglish]{babel} 
\usepackage[T1]{fontenc}

\usepackage{graphicx}
\usepackage[labelsep=endash, figurename=Figure, tablename=Table, textfont=it]{caption}
\usepackage{subcaption}
\usepackage{float}  
\usepackage{hyperref}
\usepackage[integrals]{wasysym}
\usepackage{multicol}
\usepackage{amsmath}
\usepackage{bm}

\usepackage{fancyhdr}
\usepackage[top=2.5cm, bottom=2.5cm, left=2.5cm, right=2.5cm]{geometry} 
\usepackage[usenames,dvipsnames]{pstricks}
\usepackage[round,comma,authoryear]{natbib}
\usepackage{titling}

\definecolor{darkblue}{rgb}{0,0,0.5}
\hypersetup{colorlinks=true,linkcolor=black,citecolor=darkblue,urlcolor=darkblue} 

\makeatletter
\let\ps@plain\ps@fancy
\makeatother

\begin{document}


\title{\bf Connected Vehicle Data-driven Robust Optimization for Traffic Signal Timing: Modeling Traffic Flow Variability and Errors}
\author{Chaopeng Tan$^{a}$, Yue Ding$^{b}$, Kaidi Yang$^{c}$, Hong Zhu$^{b}$, and Keshuang Tang$^{b,*}$ \\
(c.tan-2@tudelft.nl, yueding@tongji.edu.cn, kaidi.yang@nus.edu.sg, hongzhu1990@tongji.edu.cn, tang@tongji.edu.cn )}
\date{\today}

\pretitle{\centering\Large}
\posttitle{\par\vspace{1ex}}

\preauthor{\centering}
\postauthor{\par\vspace{1ex}
$^{a}$ Department of Transport and Planning, Delft University of Technology, Delft, Netherlands\\
$^{b}$ The Key Laboratory of Road and Traffic Engineering of the Ministry of Education, Tongji University, Shanghai, China\\
$^{c}$ Department of Civil and Environmental Engineering, National University of Singapore, Singapore\\
$^{*}$ Corresponding author

\vspace{1ex}\it
Extended abstract submitted for presentation at the Conference in Emerging Technologies in Transportation Systems (TRC-30)\\
September 02-03, 2024, Crete, Greece\\
\vspace{1ex}

}

\maketitle
\vspace{-1cm}
\noindent\rule{\textwidth}{0.5pt}\vspace{0cm}

\fancypagestyle{firststyle}{
\lhead[]{}
\rhead[]{}
\lfoot[TRC-30]{TRC-30}
\rfoot[Original abstract accepted]{Original abstract accepted}
\cfoot[]{}
}
\thispagestyle{firststyle}

\pagestyle{fancy}
\fancyhead{}
\fancyfoot{}
\renewcommand{\headrulewidth}{0pt}
\renewcommand{\footrulewidth}{0pt}
\setlength{\headheight}{15pt}
\lhead{Chaopeng Tan, Yue Ding, Kaidi Yang, Hong Zhu, Keshuang Tang}
\rhead[\thepage]{\thepage}
\lfoot[TRC-30]{TRC-30}
\rfoot[Original abstract accepted]{Original abstract accepted}
\cfoot[]{}


\section{ INTRODUCTION}

Traffic signal control plays an important role in improving urban traffic flow efficiency as signalized intersections stand as one of the most prevalent bottlenecks in urban road networks.
Recent advancements in Connected Vehicle (CV) technology have prompted research on leveraging CV data for more effective traffic management. 
Despite the low penetration rate, such detailed CV data has demonstrated great potential in effectively monitoring and controlling urban road networks, particularly in alleviating urban road bottlenecks. 
Specific to traffic signal control, a number of studies have been conducted to realize CV-based fixed-time control \citep{Ma2020Multi,tan2024connected} or real-time control \citep{Feng2015Real, Yang2016Isolated}. 
However, existing studies share a common shortcoming in that they all ignore traffic flow estimation errors due to the nature of CV sampling observations in their modeling process. Such errors are especially non-negligible in low penetration rate scenarios. 
To address this shortcoming, this study proposes a connected vehicle data-driven robust traffic signal optimization framework, which comprises two essential components: i) a CV-data-driven uncertainty set for arrival rates and ii) a CV-data-driven robust optimization model that explicitly accounts for traffic flow variability and potential errors. 

The \emph{major contributions} of this study are three-fold:
\underline{First}, we represent a pioneering effort in CV-based signal control research by concurrently addressing traffic flow variability and errors, which provides a paradigm for robust traffic modeling based on sampled observations, such as CVs.
\underline{Second}, we present a novel data-driven method for constructing parameter uncertainty sets using CV data. These sets, derived from the upper and lower bounds of cyclic arrival rates obtained from CV trajectories, circumvent the error-prone arrival rate estimation process.
\underline{Third}, we employ a robust optimization model to explicitly handle the uncertainties in model parameters, which can be widely implemented in different traffic scenarios, including under/over-saturated and fixed/real-time signalized intersections.


\section{METHODOLOGY}
\subsection{CV data-driven uncertainty set of arrival rate}
While it is not possible to accurately estimate the arrival rate based on CV data due to its sampling characteristics, we can infer reliable upper and lower bounds on the arrival rate based on CV trajectories during the cycle. Based on such bounds information on arrival rates, we propose a novel CV-data-driven box uncertainty set construction method.
\begin{figure}[h]
  \centering
  \includegraphics[width=0.98\textwidth]{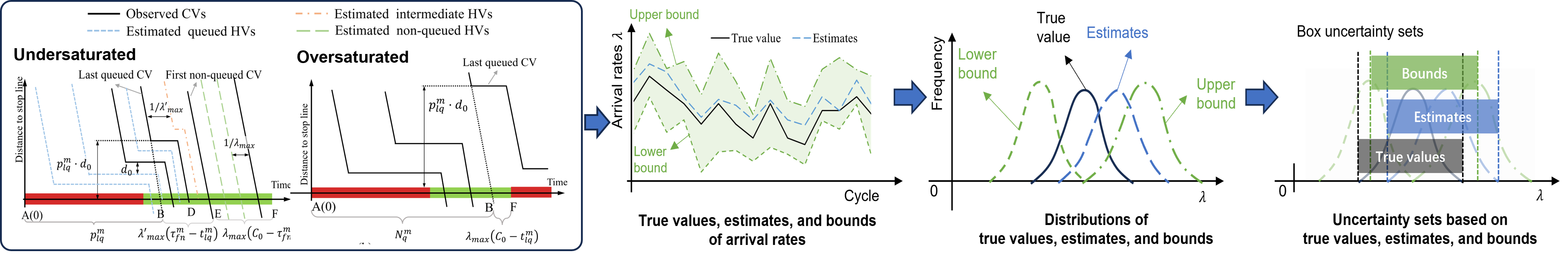}
  \caption{Illustration of the proposed CV-data-driven box uncertainty set.} \label{fig: uncertainty set}
\end{figure}
We assume that there is a set $\mathcal{M}$ of historical cycles with CV observations, i.e., those cycles of multiple days in the same time period. Note that, for fixed-time control, such a time period can be a predefined time-of-day (TOD) period, while for real-time control, such a time period can be sliding with the target cycle, e.g., a shorter five-minute period in which the cycle is located, to adapt traffic flow dynamics. As presented in Fig. \ref{fig: uncertainty set}, queued CVs provide a part of known vehicle arrivals during the cycle, which helps determine a solid lower bound $\lambda_l^m$ of the arrival rate $\lambda^m$ during the cycle. On the other hand, non-queued CVs, although do not directly provide vehicle arrival information, can help tighten the upper bound $\lambda_u^m$. Then we have,
\begin{align} 
    & \lambda_l^m = \frac{N_{1}^m + n_{nq}^m}{C_m} \leq \lambda^m \leq \lambda_u^m = \frac{N_{1}^m + \lambda'_{max}(\tau_{fn}^m - t_{lq}^m) +\lambda_{max} (C_m - \tau_{fn}^m )}{C_m}, \label{Eq: arrival rate bounds}
\end{align}
where $N_{1}^m$ is the number of arrived vehicles observed during A-B. $p$ denotes the queuing position, $t$ denotes the virtual arrival time, and $\tau$ denotes the stopline through time. The subscript $lq$, $lr$, and $fn$ denote the last queued CV, the last residual CV, and the first non-queued CV, respectively. $n_{nq}^m$ denotes the number of non-queued CVs, $C_m$ denotes the cycle length. $\lambda_{max}$ denotes the empirical maximum arrival rate and $\lambda'_{max} = \min \{\lambda_{max}, (\tau_{fn}^m - \tau_{lq}^m)/(h_s (\tau_{fn}^m - t_{lq}^m))\}$. $h_s$ is the average time headway of queued vehicles. In particular, for undersaturated cycles, we have $N_{1}^m = p_{lq}^m$. While for oversaturated cycles, we have $N_{1}^m = {(p_{lq}^m-p_{lr}^m) t_{lq}^m} /({t_{lq}^m - t_{lr}^m})$ and $\tau_{fn}^m = C_m$.
Eventually, a box uncertainty set $\mathcal{U}$ is constructed as 
\begin{align} 
    \mathcal{U} = \{ \lambda_k: \hat{l}_k \leq \lambda_k \leq \hat{u}_k, \forall k \in \mathcal{K} \}, \label{Eq: box uncertainty set}
\end{align}
where $\hat{l}_k$ and $\hat{u}_k$ are the median value of the lower and upper bounds of arrival rates for movement $k$, respectively. The reason for taking the median value rather than the mean value is that the median is less sensitive to outliers. Besides, such a box uncertainty set based on median values ensures that at least 50\% cycles are covered even in extreme cases (i.e., the arrival rates of all cycles are lower/upper bounds), which guarantees the worst coverage of the set without being too conservative as it also serves as a descriptive statistic. Since this uncertainty set is obtained from the distribution of the valid bounds, it accounts for both traffic flow variability and errors.

\subsection{CV data-driven robust signal timing optimization model}
After constructing the box uncertainty set of arrival rates based on CV data, a commonly used min-max method \citep{beyer2007robust, Yin2008Robust} is adopted for robust optimization of traffic signal timings in this study. Here we use $\Lambda$ to denote a realization of arrival rates for all movements in the proposed box uncertainty set $\mathcal{U}$; $\mathcal{F}(\theta,\Lambda)$ denotes the objective function with decision variables $\theta$ and a realization of $\Lambda$.
Then, the proposed CV-data-driven robust optimization model for traffic signal timing is written below,
\begin{align} 
    \min\limits_{\theta} \max\limits_{\Lambda \in \mathcal{U}}&  \quad \mathcal{F}(\theta,\Lambda) = \sum_{k \in \mathcal{K}} (\sum_{i \in \mathcal{I}_k} d_i + \alpha Q_k) \label{Obj robust} \\
    \text{s.t.} \quad & \text{Delay formula: }
    \begin{cases}
        d_i  \geq R_k + L_k^s - (1-\lambda_k(\Lambda)  h_k)t_i, \\
        d_i  \geq 0, \\
        i \in \mathcal{I}_k, \quad k \in \mathcal{K}, \quad \forall \Lambda \in \mathcal{U},
    \end{cases} \label{Con: CV delay, robust} \\
    & \text{Residual queue: } 
    \begin{cases}
        Q_k  \geq \lambda_k(\Lambda) C - \frac{G_k^{eff}}{h_k},\\
        Q_k  \geq 0, \\
        k \in \mathcal{K}, \quad \forall \Lambda \in \mathcal{U},
    \end{cases} \label{Con: residual queue, robust} \\
    & \text{Phase sequence: } 
    \bm{f}(\theta) \geq \bm{0}. \label{Con: phase sequence}
\end{align}
where $\theta$ denotes the set of decision variables, $\theta = \{C, \{g_k^s,g_k^e \}_{k \in \mathcal{K}} \}$, where $C$ denotes the cycle length and $g_k^s,g_k^e$ denote the green start and end times. $Q_k$ is the theoretical residual queue length. $\alpha$ is a penalty coefficient. $\alpha = T$, i.e., TOD length, for fixed-time control, and $\alpha = C_a$, i.e., historical average cycle length, for real-time control.
$\mathcal{K}$ denotes the set of controlled streams indexed by $k$, and $L_k^s, L_k^y$ denote the start-up and yellow lost times. $R_k$ denotes the red time and $R_k=C-(g_k^e-g_k^s+Y_k)$. $G_k^{eff}$ is the effective green time of movement $k$ and $G_k^{eff}=(g_k^e-g_k^s+Y_k-L_k^y-L_k^s )$. $Y_k$ denotes the yellow time.
$\mathcal{I}_k$ denotes the set of CVs in stream $k$ indexed by $i$. $d_i$ denotes its theoretical vehicle delay, $\tau_i$ is the cycle stopline through time, i.e., relative to the red start time, and $t_i$ is the cycle virtual no-delay arrival time.

In particular, for \emph{fixed-time} control, $\mathcal{I}_k$ represents historical CVs in the same TOD period, while for \emph{real-time} conrol, $\mathcal{I}_k$ refers to the CVs observed in real time on the road segment. In different cases, $t_i$ is calculated as follows after linearization, 
\begin{align}
    \text{fixed-time: }
    \begin{cases}
        t_i= t_i^0- \lfloor t_i^0/C \rfloor \cdot C - g_k^e - Y_k + b_i C, \\
        t_i^0- \lfloor t_i^0/C \rfloor \cdot C - g_k^e - Y_k + b_i M \geq 0, \\
        t_i^0- \lfloor t_i^0/C \rfloor \cdot C - g_k^e - Y_k + (b_i-1) M + \varepsilon \leq 0, \\
    \end{cases}
    \text{OR real-time: }
    t_i= t_i^0 - r_k^s
\end{align}
where $t_i^0$ is the global virtual arrival time and $r_k^s$ is the red start time. $b_i$ is a binary variable, $M$ is a big value (e.g., 300), and $\varepsilon$ is a sufficiently small value for linearization (e.g., 0.001). 

The objective Eq. (\ref{Obj robust}) minimizes the total delay of CVs, accompanied by a penalty for residual queues. Here we clarify two distinct reasons for our choice to minimize CV delays. First, when minimizing the CV delays, the proposed model can be readily linearized to facilitate the tractable solution. Second, by opting to minimize CV delays, we account for the arrival patterns exhibited by CVs to a certain extent. Constraints Eq. (\ref{Con: CV delay, robust}) derives CV delays. Constraints Eq. (\ref{Con: residual queue, robust}) is formulated to expect the intersection to be undersaturated, i.e., vehicles accumulated during the cycle can dissipate during the effective green time. Constraints Eq. (\ref{Con: phase sequence}) represent the linear restrictions imposed by the designated phase sequences, standing independent of CV data or traffic state parameters. By the robust counterpart approach, we can prove that the worst case of the proposed robust optimization is $\{\lambda_k = \Hat{u}_k\}_{k \in \mathcal{K}}$. Ultimately, the robust optimization problem can be converted to solve deterministic MILP/LP problems for either fixed/real-time control.

\section{REPRESENTATIVE RESULTS}
Due to page limitations, here we only present representative results for the proposed method at a real-world fixed-time signalized intersection situated at Yuqiao Middle Road-Qingfeng North Road in Tongxiang City, in China. Based on the traffic volume data collected in the field, we calibrated a SUMO micro-simulation model to test the effectiveness of the different methods, i.e., Synchro, CV-based deterministic optimization model (CV-DO, use average arrival rate estimates as input), and the proposed CV-based robust optimization model (CV-RO). 
\begin{figure}[h]
  \centering
  \includegraphics[width=0.77\textwidth]{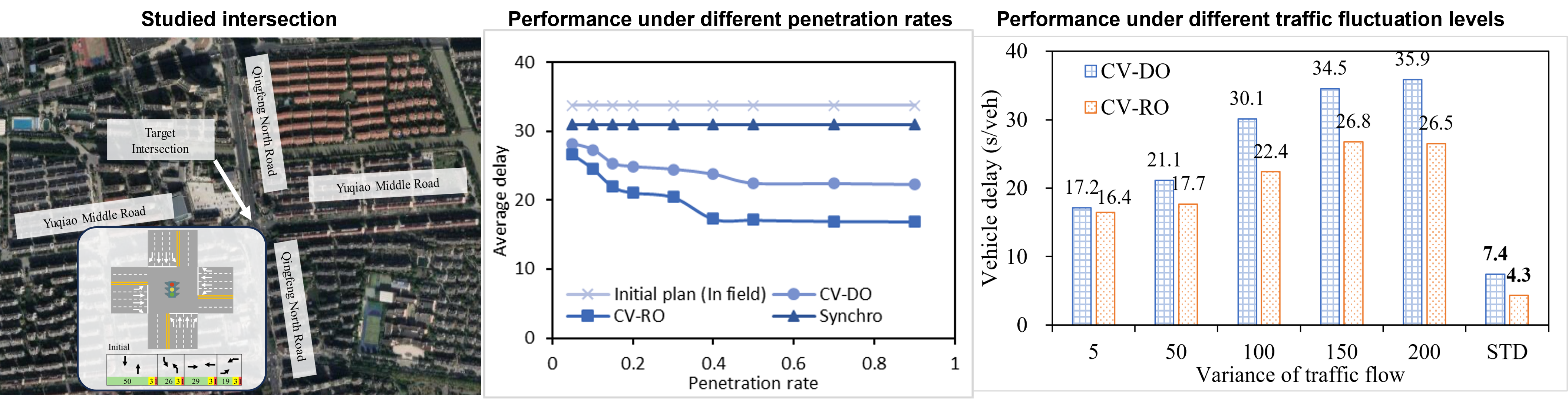}
  \caption{Representative test results at a fixed-time signalized intersection.} \label{fig: results}
\end{figure}
Fig. \ref{fig: results} presents the performance of benchmark methods under various penetration rates and traffic flow fluctuation levels. (a) As the penetration rate increases, both CV-based optimization models demonstrate a consistent decrease in average vehicle delay. (b) When compared to Synchro, which relies on precise arrival rates as input, the CV-RO yields lower average vehicle delay at all penetration rates. (c) CV-RO outperforms CV-DO significantly at higher penetration rates. (d) CV-RO consistently achieves lower vehicle delays compared to CV-DO across different traffic flow fluctuation levels and such advantages are more pronounced in scenarios characterized by greater variability.

\section{DISCUSSION}
We proposed a CV-based robust optimization (CV-RO) method for fixed/real-time signalized intersections, aimed at addressing the challenges posed by traffic flow variability and inherent estimation errors using CV data. Evaluation results highlight the superior performance of the CV-RO model compared to benchmark models across various CV penetration rates and traffic flow fluctuation levels, affirming its effectiveness and robustness. In particular, the study provides \emph{inspiration for future research on how to proactively address potential parameter errors from incomplete observations for traffic models}, as suggested by the proposed method that is entirely data-driven and does not require pre-calibration of the parameter error distributions through validation data.


\begin{small}
\begin{sloppypar} 
\bibliographystyle{authordate1} 

\setlength{\bibsep}{0pt}

\bibliography{Manuscript_V1}

\end{sloppypar}
\end{small}


\end{document}